\documentclass{article}
\usepackage{amsmath,amsthm,amsopn,amstext,amscd,amsfonts,amssymb}
\usepackage{natbib}
\usepackage{verbatim}

\def\diag{\mathop{\rm diag}\nolimits}
\def\vol{\mathop{\rm Vol}\nolimits}
\def\tr{\mathop{\rm tr}\nolimits}
\def\R{\mathop{\rm Re}\nolimits}

\def\etr{\mathop{\rm etr}\nolimits}

\renewenvironment{abstract}
                 {\vspace{6pt}
                  \begin{center}
                  \begin{minipage}{5in}
                  \centerline{\textbf{Abstract}}
                  \noindent\ignorespaces
                 }
                 {\end{minipage}\end{center}}
\newtheorem{thm}{\textbf{Theorem}}[section]
\newtheorem{cor}{\textbf{Corollary}}[section]
\newtheorem{lem}{\textbf{Lemma}}[section]
\theoremstyle{definition}
\newtheorem{defn}{\textbf{Definition}}[section]

\newtheorem{rem}{\textbf{Remark}}[section]

\title{\huge \textbf{Complex bimatrix variate generalised beta distributions}}
\author{
  \textbf{Jos\'e A. D\'{\i}az-Garc\'{\i}a} \thanks{Corresponding author\newline
    }\\
  Department of Statistics and Computation \\
  25350 Buenavista, Saltillo, Coahuila, Mexico \\
  E-mail: jadiaz@uaaan.mx \\[2ex]
  \textbf{Ram\'on Guti\'errez J\'aimez} \\
  Department of Statistics and O.R. \\
  University of Granada \\
  Granada 18071, Spain \\
  E-mail: rgjaimez@ugr.es\\
}
\date{}
\begin{document}
\maketitle

\begin{abstract}
In this paper, the study of bivariate generalised beta type I and II distributions is
extended to the complex matrix variate case, for which the corresponding density
functions are found. In addition, for complex bimatrix variate beta type I
distributions, several basic properties, including the joint eigenvalue density and
the maximum eigenvalue distribution, are studied.
\end{abstract}

\section{Introduction}

Complex matrix variate distributions play an important role in various fields of
research. Applications of complex random matrices can be found in multiple time
series analysis, nuclear physics, complex multivariate linear model, shape theory,
the evaluation of the capacity of multiple-input multiple-output (MIMO) wireless
communication systems, see \citet{m:91}, \citet{k:65}, \citet{mdm:06} and
\citet{rva:05a} among many others.

Several studies on the distribution of complex random matrices have been made, see
\citet{mp:05}. The complex matrix variate Gaussian distribution was introduced by
\citet{w:56}, and further developed by \citet{t:60} and \citet{g:63}. The complex
Wishart distribution was studied by \citet{g:63} and \citet{j:64}, among many others.
\citet{j:64} and \citet{k:65} derived  both complex central and noncentral matrix
variate beta distributions.

When the goal is to generalise the distribution of a random variable to the
multivariate case, two options are normally addressed, by which it is extended to
either the vectorial or the matrix case, e.g. normal, t or bessel distributions,
among many others. However, some of these generalisations have traditionally been
made directly to the matrix case, when such a matrix is symmetric, as in the case of
the chi-square and beta distributions, for which the corresponding multivariate
distributions are the Wishart and matrix variate beta distributions, respectively.
Nevertheless, these latter generalisations are inappropriate in some cases, because,
in some applications, the researcher is interested in a vector variate, not in a
symmetric matrix, see \citet{ln:82}. In other words, the researcher is interested in
a vector, say, $\mathbf{X} = (x_{1}, \dots , x_{m})'$, such that $x_{i}$ has a
marginal beta type I or II distribution for all $i = 1, \dots,m$. In this respect,
\citet{ln:82} and \citet{cn:84} proposed a multivariate vector-version of the beta
type I and II distributions. Let us consider the following bivariate version, see
\citet{ol:03}.

Let $X_{0}, X_{1},X_{2}$ be distributed as independent gamma random variates with
parameters $a = a_{0}, a_{1}, a_{2}$, respectively (see Definition \ref{defgama} in
Section 2); and define
\begin{equation}\label{eq1}
    U_{1} = \frac{X_{1}}{X_{1}+X_{0}}, \qquad U_{2} = \frac{X_{2}}{X_{2}+X_{0}}.
\end{equation}
Clearly, $U_{1}$ and $U_{2}$ each have a beta type I distribution, $U_{1} \sim
\mathcal{B}I_{1}(a_{1},a_{0})$ and $U_{2} \sim \mathcal{B}I_{1}(a_{2},a_{0})$, over
$0 \leq u_{1},u_{2} \leq 1$ (see Subsection 2.1). However, they are correlated such
that $(U_{1},U_{2})'$ has a bivariate generalised beta type I distribution over $0
\leq u_{1},u_{2} \leq 1$. The kernel of the joint density function of $U_{1}$ and
$U_{2}$ is
$$
  \propto \frac{u_{1}^{a_{1}-1} u_{2}^{a_{2}-1} (1 - u_{1})^{a_{2} + a_{0}-1} (1 - u_{2})^{a_{1} + a_{0}-1}}
  {(1 - u_{1}u_{2})^{a_{1}+ a_{2}+ a_{0}}}, \quad 0 \leq u_{1},u_{2} \leq 1.
$$

A similar result is obtained in the case of beta type II. Here it defines
$$
  F_{1} = \frac{X_{1}}{X_{0}}, \qquad F_{2} = \frac{X_{2}}{X_{0}}.
$$
Once again it is evident that $F_{1}$ and $F_{2}$ each have a beta type II
distribution, $F_{1} \sim \mathcal{B}II_{1}(a_{1},a_{0})$ and $F_{2} \sim
\mathcal{B}II_{1}(a_{2},a_{0})$, over $f_{1},f_{2} \geq 0$. As in the beta type I
case, they are correlated such that $(F_{1},F_{2})'$ has a bivariate generalised beta
type II distribution over $f_{1},f_{2} \geq 0$. The kernel of the joint density
function of $F_{1}$ and $F_{2}$ is
$$
  \propto \frac{f_{1}^{a_{1}-1} f_{2}^{a_{2}-1}}
  {(1 + f_{1} + f_{2})^{a_{1}+ a_{2}+ a_{0}}},
  \quad f_{1},f_{2} \geq 0.
$$
Some applications to utility modelling  and Bayesian analysis are presented in
\citet{ln:82} and \citet{cn:84}, respectively.

These ideas can be extended to the matrix variate case. Thus, let us assume a
partitioned matrix $\mathbb{U} = (\mathbf{U}_{1} \vdots \mathbf{U}_{2})' \in
\mathfrak{C}^{2m \times m}$, then under the complex matrix variate versions of the
transformations (\ref{eq1}), we are interested in finding the joint density of
$\mathbf{U}_{1}$ and $\mathbf{U}_{2}$, from where it is easy to see that the marginal
densities of $\mathbf{U}_{1}$ and $\mathbf{U}_{2}$ are complex matrix variate beta
type I distributions. In the central and noncentral real cases, the matrix variate
joint densities of $\mathbf{U}_{1}$ and $\mathbf{U}_{2}$ and of $\mathbf{F}_{1}$ and
$\mathbf{F}_{2}$, together with some of their properties, are studied in
\citet{dggj:09a, dggj:09b}. These distributions are termed central complex bimatrix
variate generalised beta type I and II distributions, respectively.

In the present paper, the bivariate generalised beta type I and II distributions are
extended the the complex matrix variate case, see Sections \ref{sec3} and \ref{sec4}.
In Section \ref{sec5}, certain basic properties, the joint eigenvalue density and the
density of their maximum eigenvalues are studied for the complex bimatrix variate
generalised beta type I distribution.

\section{Preliminary results}\label{sec2}

The join eigenvalues density can be calculated using complex hypergeometric functions
and complex invariant polynomials with matrix arguments.  In this section, these and
other issues are addressed for the case of the complex multivariate distribution. Let
us first establish some notation.

\subsection{Notation and matrix variate distributions}

Let $\mathbf{A}= (a_{rs}) = (a_{rs_{1}} + i a_{rs_{2}}) = \mathbf{A}_{1} + i
\mathbf{A}_{2}$ be an $m \times n$ of complex numbers, $\mathbf{A} \in
\mathfrak{C}^{m \times n}$, where $\mathbf{A}_{1}= (a_{rs_{1}})$ and $\mathbf{A}_{2}=
(a_{rs_{2}})$ are real matrices $m \times n$, $\mathbf{A}_{1}, \mathbf{A}_{2} \in
\Re^{m \times n}$, and $i = \sqrt{-1}$. Then, $\mathbf{A}'$ denotes the transpose of
$\mathbf{A}$, $\overline{\mathbf{A}}$ denotes the conjugate of $\mathbf{A}$, and
$\mathbf{A}^{H}$ denotes the conjugate transpose of $\mathbf{A}$. For $n = m$, let
$\tr(\mathbf{A}) = a_{11} + \cdots +a_{mm}$, $\etr(\tr(\mathbf{A})) =
\exp(\tr(\mathbf{A}))$, and then $|\mathbf{A}|$ denotes the determinant of
$\mathbf{A}$. $\mathbf{A} = \mathbf{A}^{H} > \mathbf{0}$ is a Hermitian positive
definite matrix, and $\mathbf{A}^{1/2}$ denotes the unique Hermitian positive
definite square root matrix of $\mathbf{A} = \mathbf{A}^{H} > \mathbf{0}$.
$d\mathbf{A}$ denotes the differential matrix of $\mathbf{A}$ and $(d\mathbf{A}) =
(d\mathbf{A}_{1})(d\mathbf{A}_{2})$ denotes the volume element (Lebesgue measure)
commonly associated with $\mathbf{A}$. For example if $\mathbf{A}$ is a Hermitian
matrix, then $\mathbf{A}_{1}$ is a symmetric matrix ($\mathbf{A} = \mathbf{A}'$) and
$\mathbf{A}_{2}$ is a skew-symmetric matrix ($\mathbf{A} = -\mathbf{A}'$). In this
case
$$
  (d\mathbf{A}) = \bigwedge_{r \leq s}da_{rs_{1}} \bigwedge_{r < s}
  da_{rs_{2}}.
$$

The space of all matrices $\mathbf{G}_{1} \in \mathfrak{C}^{n \times m}$ ($m \leq n$)
with orthonormal columns is termed the \textbf{Stiefel manifold}, denoted by
$\mathfrak{C}\mathcal{V}_{m,n}$. Thus
$$
  \mathfrak{C}\mathcal{V}_{m,n} = \{\mathbf{G}_{1} \in \mathfrak{C}^{n \times m}|
  \mathbf{G}_{1}^{H} \mathbf{G}_{1} = \mathbf{I}_{m}\}.
$$
From \citet{j:64}
$$
  \vol(\mathfrak{C}\mathcal{V}_{m,n}) = \int_{\mathbf{G}_{1} \in \mathfrak{C}\mathcal{V}_{m,n}}
  (\mathbf{G}_{1}^{H}d\mathbf{G}_{1}) = \frac{2^{m}
  \pi^{mn}}{\mathfrak{C}\Gamma_{m}[n]},
$$
where $\mathfrak{C}\Gamma_{m}[a]$ denotes the complex multivariate gamma function and
is defined as
$$
  \mathfrak{C}\Gamma_{m}[a] = \int_{\mathbf{V}=\mathbf{V}^{H}> \mathbf{0}} \etr(-\mathbf{V}) |\mathbf{V}|^{a-m} (d\mathbf{V})
  = \pi^{m(m-1)/2} \prod_{j = 1}^{m} \Gamma[a -j+1],
$$
where $\R(a) > m-1$.

If $m = n$, it is have a special case of the Stiefel manifold, termed unitary
manifold or unitary group and denoted as $\mathcal{U}(m) \equiv
\mathfrak{C}\mathcal{V}_{m,m}$.

\begin{defn}[The complex multivariate beta function]\label{defgama}
The complex multivariate beta function, denoted as
$\mathfrak{C}\boldsymbol{\beta}_{m}[a,b]$, is defined by
\begin{eqnarray*}
  \mathfrak{C}\boldsymbol{\beta}_{m}[b,a] &=& \int_{\mathbf{0}<\mathbf{S}=\mathbf{S}^{H}<\mathbf{I}_{m}} |\mathbf{S}|^{a-m}
  |\mathbf{I}_{m} - \mathbf{S}|^{b-m} (d\mathbf{S}) \\
    &=& \int_{\mathbf{R}=\mathbf{R}^{H}>\mathbf{0}} |\mathbf{R}|^{a-m} |\mathbf{I}_{m} + \mathbf{R}|^{-(a+b)} (d\mathbf{R}), \
    \mathbf{R} = (\mathbf{I}-\mathbf{S})^{-1}-\mathbf{I} \\
    &=& \frac{\mathfrak{C}\Gamma_{m}[a] \mathfrak{C}\Gamma_{m}[b]}{\mathfrak{C}\Gamma_{m}[a+b]}.
\end{eqnarray*}
where $\R(a) > m-1$ and $\R(b) > m-1$.
\end{defn}
We now give definitions of the complex matrix variate gamma, beta type I and II
distributions, see \citet{j:64}, \citet{k:65} and \citet{m:97}.

\begin{defn}[Complex matrix variate gamma distribution] It is said that $\mathbf{A} \in \mathfrak{C}^{m \times
m}$,  a random Hermitian positive definite, has a complex matrix variate gamma
distribution with parameters $a$ and a Hermitian positive definite matrix
$\mathbf{\Theta} \in \mathfrak{C}^{m \times m}$, if its density function is
\begin{equation}\label{gamma}
    \frac{1}{\mathfrak{C}\Gamma_{m}[a]|\mathbf{\Theta}|^{a}}|\mathbf{A}|^{a-m}\etr(-\mathbf{\Theta}^{-1}\mathbf{A})
    (d\mathbf{A}), \quad \mathbf{A}=\mathbf{A}^{H} > \mathbf{0},
\end{equation}
where $\R(a) > m-1$. Such a distribution is  denoted as $\mathbf{A} \sim
\mathfrak{C}\mathcal{G}_{m}(a, \mathbf{\Theta})$.
\end{defn}

\begin{lem}[Complex matrix variate beta type I and II distribution]
If $\mathbf{A}$ and $\mathbf{B}$ have a complex matrix variate gamma distribution,
i.e. $\mathbf{A} \sim \mathfrak{C}\mathcal{G}_{m}(a,\mathbf{I}_{m})$ and $\mathbf{B}
\sim \mathfrak{C}\mathcal{G}_{m}(b,\mathbf{I}_{m})$ independently.
\begin{enumerate}
\item Then the complex matrix variate beta type I distribution is defined as
\begin{equation}\label{defbI}
  \mathbf{U} =
  \left\{%
     \begin{array}{ll}
      (\mathbf{A} + \mathbf{B})^{-1/2}\mathbf{A} ((\mathbf{A} + \mathbf{B})^{-1/2})', & \mbox{Definition $1$ or},\\
       \mathbf{A} ^{1/2}(\mathbf{A} + \mathbf{B})^{-1} (\mathbf{A} ^{1/2})', & \mbox{Definition $2$}.\\
    \end{array}%
  \right.
\end{equation}
Thus under definitions $1$ and $2$ its density function is denoted as
$$
  \mathfrak{C}\mathcal{B}I_{m}(\mathbf{U};a,b),
$$
and given by
\begin{equation}\label{beta}
     \frac{1}{\mathfrak{C}\boldsymbol{\beta}_{m}[a,b]} |\mathbf{U}|^{a-m} |\mathbf{I}_{m} - \mathbf{U}|^{b - m}
     (d\mathbf{U}), \quad \mathbf{0} < \mathbf{U}= \mathbf{U}^{H} < \mathbf{I}_{m},
\end{equation}
this being denoted as $\mathbf{U} \sim \mathfrak{C}\mathcal{B}I_{m}(a,b)$ with
Re$(a)> m-1$ and Re$(b)> m-1$.
\item Then the complex matrix variate beta type II distribution is defined as
\begin{equation}\label{defbII}
  \mathbf{F} =
  \left\{%
     \begin{array}{ll}
      \mathbf{B}^{-1/2}\mathbf{A} (\mathbf{B}^{-1/2})', & \mbox{Definition 1},\\
      \mathbf{A}^{1/2}\mathbf{B}^{-1} (\mathbf{A} ^{1/2})', & \mbox{Definition 2},\\
    \end{array}%
  \right.
\end{equation}
Thus under definitions $1$ and $2$ its density function is denoted as
$$
  \mathfrak{C}\mathcal{B}II_{m}(\mathbf{U};a,b),
$$
and given by
\begin{equation}\label{efe}
    \frac{1}{\mathfrak{C}\boldsymbol{\beta}_{m}[a,b]}  |\mathbf{F}|^{a-m}|\mathbf{I}_{m} +
\frac{}{}    \mathbf{F}|^{-(a+b)}.
    (d\mathbf{F}), \quad \mathbf{F} = \mathbf{F}^{H} > \mathbf{0}.
\end{equation}
this being denoted as $\mathbf{F} \sim \mathfrak{C}\mathcal{B}II_{m}(a,b)$ with
Re$(a)> m-1$ and Re$(b)> m-1$.
\end{enumerate}
\end{lem}

\subsection{Complex hypergeometric functions and invariant polynomials}

The following definitions of hypergeometric functions with a matrix argument are
based on \citet{c:63} and \citet{ke:06}.

\begin{defn}\label{def1}
The hypergeometric functions of a matrix argument are given by
\begin{equation}\label{hf}
    {}_{p}F_{q}^{(\alpha)}(a_{1},\dots,a_{p};b_{1},\dots,b_{q};\mathbf{X}) = \sum_{t=0}^{\infty} \sum_{\tau}
    \frac{(a_{1})^{(\alpha)}_{\tau} \cdots (a_{p})^{(\alpha)}_{\tau}}{(b_{1})^{(\alpha)}_{\tau} \cdots (b_{q})^{(\alpha)}_{\tau}}
    \frac{C_{\tau}^{(\alpha)}(\mathbf{X})}{t!},
\end{equation}
where $\sum_{\tau}$ denotes the summation over all the partitions $\tau = (t_{1},
\dots, t_{m})$, $t_{1} \geq \cdots \geq t_{m} \geq 0$, of $t=t_{1}+\cdots + t_{m}$,
$C^{(\alpha)}_{\tau}(\mathbf{X})$ is the Jack polynomial of $\mathbf{X}$
corresponding to $\tau$ and the generalised hypergeometric coefficient
$(a)^{(\alpha)}_{\tau}$ is given by
$$
  (a)^{(\alpha)}_{\tau} = \prod_{j = 1}^{m} (a - (i-1)/(\alpha))_{t_{i}},
$$
where $(a)_{t} = a(a+1)(a+2) \cdots (a+t-1)$, $(a)_{0} = 1$. Here $\mathbf{X} \in
\mathfrak{C}^{m \times m}$, is a complex symmetric matrix and the parameters $a_{i}$,
$b_{j}$ are arbitrary complex numbers.
\end{defn}
Other characteristics of the parameters $a_{i}$ and $b_{j}$ and the convergence of
(\ref{hf}) appear in \citet[p. 258]{mh:82}, \citet{gr:87} and \citet{rva:05b}.

\begin{rem}
In Definition \ref{def1}, when $\alpha =1$ and $2$ the complex and real cases are
obtained, respectively. In this paper it is considered only the complex case. Then,
adopting the notation used by \citet{j:64}, denoting the complex hypergeometric
function, generalised hypergeometric coefficient and zonal polynomials as
${}_{p}\widetilde{F}_{q}^{(\alpha)} \equiv {}_{p}F_{q}^{(1)}$, $[a]_{\tau} \equiv
(a)_{\tau}^{(1)}$, and $\widetilde{C}_{\tau}(\cdot) \equiv C_{\tau}^{(1)}(\cdot)$,
respectively
\end{rem}

A special case of (\ref{hf}) is
\begin{eqnarray*}
  {}_{1}\widetilde{F}_{0}(a;\mathbf{X}) &=& \sum_{t=0}^{\infty} \sum_{\tau}
    [a]_{\tau}  \frac{\widetilde{C}_{\tau}(\mathbf{X})}{t!}, \quad (||\mathbf{X}|| < 1)\\
    &=& |\mathbf{I}_{m} - \mathbf{X}|^{-a},
\end{eqnarray*}
where $\|\mathbf{X}\|$ denotes the maximum of the absolute values of the eigenvalues
of $\mathbf{X}$.

From \citet{rva:05b}, is known that,
\begin{eqnarray*}
    {}_{1}\widetilde{F}_{1}(a;c; \mathbf{X})=
    \frac{1}{\mathfrak{C}\boldsymbol{\beta}_{m}[a,c-a]}\hspace{6cm}\\
    \times \int_{\mathbf{0}<\mathbf{Y} = \mathbf{Y}^{H}<\mathbf{I}_{m}} {}_{0}\widetilde{F}_{0}(\mathbf{XY})|\mathbf{Y}|^{a-m}
    |\mathbf{I}-\mathbf{Y}|^{c-a-m}(d\mathbf{Y}),
\end{eqnarray*}
and
\begin{eqnarray}
    {}_{2}\widetilde{F}_{1}(a,a_{1};c; \mathbf{X})=
    \frac{1}{\mathfrak{C}\boldsymbol{\beta}_{m}[a,c-a]}\hspace{6cm}\nonumber\\ \label{euler}
    \times \int_{\mathbf{0}<\mathbf{Y} = \mathbf{Y}^{H}<\mathbf{I}_{m}} {}_{1}\widetilde{F}_{0}(a_{1};\mathbf{XY})|\mathbf{Y}|^{a-m}
    |\mathbf{I}-\mathbf{Y}|^{c-a-m}(d\mathbf{Y}).
\end{eqnarray}
Thus we have
\begin{lem}\label{lemma}
Let $\mathbf{X}= \mathbf{X}^{H} \in \mathfrak{C}^{m \times m}$, with
$\R(\mathbf{X})<\mathbf{I}$, $\R(a)>m-1$, $\R(c)>m-1$ and $\R(c-a)>m-1$. Then
\begin{eqnarray*}
    {}_{p+1}\widetilde{F}_{q+1}(a,a_{1}, \dots,a_{p};c,b_{1}, \dots, b_{q}; \mathbf{X})=
    \frac{1}{\mathfrak{C}\boldsymbol{\beta}_{m}[a,c-a]}\hspace{4cm}\\
    \times \int_{\mathbf{0}<\mathbf{Y} = \mathbf{Y}^{H}<\mathbf{I}_{m}} {}_{p}\widetilde{F}_{q}(a_{1} \cdots a_{p};b_{1} \cdots b_{q};
    \mathbf{XY})|\mathbf{Y}|^{a-m} |\mathbf{I}-\mathbf{Y}|^{c-a-m}(d\mathbf{Y}).
\end{eqnarray*}
\end{lem}
\proof First, an expansion is applied in terms of complex zonal polynomials
$$
    {}_{p}\widetilde{F}_{q}(a_{1}, \dots,a_{p};b_{1}, \dots, b_{q}; \mathbf{XY})=\sum_{t=0}^{\infty}
    \sum_{\tau} \frac{[a_{1}]_{\tau} \cdots[a_{p}]_{\tau}}{[b_{1}]_{\tau} \cdots [b_{q}]_{\tau}} \frac{\widetilde{C}_{\tau}(\mathbf{XY})}{t!}.
$$
Then, after integrating term by term, see \citet{rva:05b}, we have that
\begin{eqnarray*}
    &&\hspace{-1cm}\int_{\mathbf{0}<\mathbf{Y} = \mathbf{Y}^{H}<\mathbf{I}_{m}}{}_{p}\widetilde{F}_{q}(a_{1} \cdots a_{p};b_{1} \cdots b_{q};
    \mathbf{XY})|\mathbf{Y}|^{a-m}|\mathbf{I}- \mathbf{Y}|^{c-a-m}(d\mathbf{Y}) \\
    &&\hspace{-0.7cm}=\sum_{t=0}^{\infty}
    \sum_{\tau} \frac{[a_{1}]_{\tau} \cdots [a_{p}]_{\tau}}{[b_{1}]_{\tau} \cdots [b_{q}]_{\tau} \ t!}
    \int_{\mathbf{0}<\mathbf{Y} = \mathbf{Y}^{H}<\mathbf{I}_{m}}|\mathbf{Y}|^{a-m}|\mathbf{I}-\mathbf{Y}|^{c-a-m}
    \widetilde{C}_{\tau}(\mathbf{XY})(d\mathbf{Y})\\
    &&\hspace{-0.7cm}=  \mathfrak{C}\boldsymbol{\beta}_{m}[a,c-a]\sum_{t=0}^{\infty}\sum_{\tau} \frac{[a]_{\tau}}{[c]_{\tau}}
    \frac{[a_{1}]_{\tau} \dots [a_{p}]_{\tau}}{[b_{1}]_{\tau} \cdots [b_{q}]_{\tau} \ t!}
    \widetilde{C}_{\tau}(\mathbf{X})\\
    &&\hspace{-0.7cm}= \mathfrak{C}\boldsymbol{\beta}_{m}[a,c-a]\,  {}_{p+1}\widetilde{F}_{q+1}(a,a_{1} \cdots a_{p};c,b_{1}
    \cdots b_{q}; \mathbf{X}),
\end{eqnarray*}
and the required result follows. \qed

The use of complex zonal polynomials and the hypergeometric function with a matrix
argument has only recently been extended; to a large extent this is derived from the
work of \citet{ke:06}, who in \citet{k:04} provided a program in MatLab with a very
efficient algorithm for calculating Jack polynomials (in particular complex zonal
polynomials) and the complex hypergeometric function with a matrix argument.

This section concludes by establishing the following two properties of a class of
homogeneous polynomials $\widetilde{C}^{\kappa,\tau}_{\phi}(\mathbf{R},\mathbf{S})$
of degrees $k$ and $t$ in the eigenvalues of the Hermitian matrices $\mathbf{R},
\mathbf{S} \in \mathfrak{C}^{m \times m}$, respectively, see \citet{da:80} and
\citet{rva:05a}. These properties generalise the incomplete beta function equation
(61) of \citet{c:63}. The first is proposed by \citet[eq. (3.3)]{da:79} and the
second is obtained using the complex versions of \citet[eq. (3.33)]{ch:80} and the
review version of \citet[eq. (3.11)]{ch:80} given in \citet[eq. (2.7)]{chd:86}.

\begin{lem}\label{lem3}
Let $\mathbf{R}, \mathbf{S}, \mathbf{\Omega}, \mathbf{\Xi} \in \mathfrak{C}^{m \times
m}$ Hermitian matrices. Then
$$
  \int_{\mathbf{0}}^{\mathbf{\Delta}=\mathbf{\Delta^{H}}} |\mathbf{R}|^{a-m} |\mathbf{I} -
  \mathbf{R}|^{b-m} \widetilde{C}_{\tau}(\mathbf{\Omega}\mathbf{R})
  (d\mathbf{R})\hspace{6cm}
$$
\begin{equation}\label{ib1}
  =  \mathfrak{C}\boldsymbol{\beta}_{m}[a,m]|\mathbf{\Delta}|^{a}\sum_{k=0}^{\infty}
  \ \sum_{\kappa;\phi \in \kappa.\tau}\frac{[a]_{\phi}[-b+m]_{\kappa} \theta_{\phi}^{\kappa, \tau}
  \widetilde{C}_{\phi}^{\kappa, \tau}(\mathbf{\Delta}, \mathbf{\Omega}\mathbf{\Delta})}{k! [a+m]_{\phi}}
\end{equation}
and
$$
  \int_{\mathbf{0}}^{\mathbf{\Delta}=\mathbf{\Delta^{H}}} |\mathbf{R}|^{a-m} |\mathbf{I} -
  \mathbf{R}|^{b-m} \widetilde{C}_{\phi}^{\kappa, \tau}(\mathbf{\Xi},\mathbf{\Omega}\mathbf{R})
  (d\mathbf{R})\hspace{6cm}
$$
\begin{equation}\label{ib2}
  =  \mathfrak{C}\boldsymbol{\beta}_{m}[a,m]|\mathbf{\Delta}|^{a}\sum_{s=0}^{\infty}
  \ \sum_{\sigma;\phi_{1} \in \sigma. \phi}\frac{[a]_{\phi_{1}}[-b+m]_{\sigma}
  \pi_{\sigma,\phi}^{\sigma, \kappa, \tau;\phi_{1}}
  \widetilde{C}_{\phi_{1}}^{\sigma, \kappa, \tau}(\mathbf{\Delta}, \mathbf{\Xi},
  \mathbf{\Omega}\mathbf{\Delta})}{s! [a+m]_{\phi_{1}}}
\end{equation}
where, from \citet[Lemma 2.2(i) and (ii), respectively]{chd:86} it is known that
$$
  \theta_{\phi}^{\kappa, \tau} = \frac{\widetilde{C}_{\phi}^{\kappa,
  \tau}(\mathbf{I},\mathbf{I})}{\widetilde{C}_{\phi}(\mathbf{I})}\quad \mbox{and}\quad
  \pi_{\sigma,\phi}^{\sigma, \kappa, \tau;\phi_{1}} = \sum_{\phi'_{1} \equiv \phi_{1}}
  \gamma_{\sigma;\phi_{1}'}^{\sigma, \kappa, \tau;\phi_{1}}
  \overline{\alpha}_{\phi}^{\sigma^{*},\kappa, \tau;\phi'_{1}}.
$$
\end{lem}

\section{Bimatrix variate generalised beta type I distribution}\label{sec3}

Let $\mathbf{A}$, $\mathbf{B}$ and $\mathbf{C}$ be independent complex random
matrices, such that $\mathbf{A} \sim \mathfrak{C}\mathcal{G}_{m}(a, \mathbf{I}_{m})$,
$\mathbf{B} \sim \mathfrak{C}\mathcal{G}_{m}(b, \mathbf{I}_{m})$ and $\mathbf{C} \sim
\mathfrak{C}\mathcal{G}_{m}(c, \mathbf{I}_{m})$ with $\R(a)
> m-1$, $\R(b) > m-1$ and $\R(c) > m-1$ and let us define
\begin{equation}\label{bgb1}
    \mathbf{U}_{1} = (\mathbf{A}+\mathbf{C})^{-1/2}\mathbf{A}(\mathbf{A} + \mathbf{C})^{-1/2}
  \quad \mbox{and} \quad \mathbf{U}_{2} = (\mathbf{B}+\mathbf{C})^{-1/2}\mathbf{B}
  (\mathbf{B}+\mathbf{C})^{-1/2}.
\end{equation}
Of course, $\mathbf{U}_{1} \sim \mathfrak{C}\mathcal{B}I_{m}(a,c)$ and
$\mathbf{U}_{2} \sim \mathfrak{C}\mathcal{B}I_{m}(b,c)$. However, they are correlated
such that the distribution of $\mathbb{U} = (\mathbf{U}_{1}\vdots \mathbf{U}_{2})'
\in \mathfrak{C}^{2m \times m}$ can be termed a complex bimatrix variate generalised
beta type I distribution, denoted as
$$
  \mathbb{U} \sim \mathfrak{C}\mathcal{BGB}I_{2m \times m}(a,b,c).
$$

\begin{thm}\label{teob1}
Assume that $\mathbb{U} \sim \mathcal{BGB}I_{2m \times m}(a,b,c)$. Then its density
function is
\begin{equation}\label{density1}\hspace{-0.7cm}
  \frac{|\mathbf{U}_{1}|^{a-m} |\mathbf{U}_{2}|^{b-m} |\mathbf{I}_{m} - \mathbf{U}_{1}|
  |^{b + c -m} |\mathbf{I}_{m} - \mathbf{U}_{2}|^{a + c -m}}{\mathfrak{C}\boldsymbol{\beta}^{*}_{m}[a,b,c]
  |\mathbf{I}_{m} - \mathbf{U}_{1}\mathbf{U}_{2}|^{a + b + c }}(d\mathbb{U}),
\end{equation}
$\mathbf{0} < \mathbf{U}_{1} = \mathbf{U}_{1}^{H}< \mathbf{I}_{m}$, $\mathbf{0} <
\mathbf{U}_{2} = \mathbf{U}_{2}^{H}< \mathbf{I}_{m}$, where the measure
$$
  (d\mathbb{U}) = (d\mathbf{U}_{1})\wedge (d\mathbf{U}_{2}),
$$
and
$$
  \mathfrak{C}\boldsymbol{\beta}^{*}_{m}[a,b,c] = \frac{\mathfrak{C}\Gamma_{m}[a] \mathfrak{C}\Gamma_{m}[b]
  \mathfrak{C}\Gamma_{m}[c]}{\mathfrak{C}\Gamma_{m}[a+b+c]}.
$$
and $\R(a)> m-1$, $\R(b) > m-1$ and $\R(c) > m-1$.
\end{thm}
\proof The joint density of $\mathbf{A}$, $\mathbf{B}$ and $\mathbf{C}$ is
$$
  \frac{|\mathbf{A}|^{a-m} |\mathbf{B}|^{b-m} |\mathbf{C}|^{c-m}}{\mathfrak{C}\Gamma_{m}[a]
  \mathfrak{C}\Gamma_{m}[b] \mathfrak{C}\Gamma_{m}[c]} \etr(-(\mathbf{A} + \mathbf{B} + \mathbf{C})) (d\mathbf{A})
  (d\mathbf{B})(d\mathbf{C}).
$$
By effecting the change of variable (\ref{bgb1}), and taking into account
\citet[Theorems 3.5 and 3.8, pp. 183 and 190, respectively]{m:97} we have
$$
  (d\mathbf{A})(d\mathbf{B})(d\mathbf{C}) = |\mathbf{C}|^{2m} |\mathbf{I}_{m} - \mathbf{U}_{1}|
  ^{-2m} |\mathbf{I}_{m} - \mathbf{U}_{2}|^{-2m}
  (d\mathbf{U}_{1})(d\mathbf{U}_{2})(d\mathbf{C}).
$$
The joint density of $\mathbf{U}_{1}$, $\mathbf{U}_{2}$ and $\mathbf{C}$ is
$$
  \frac{|\mathbf{U}_{1}|^{a-m} |\mathbf{U}_{2}|^{b-m} }{\mathfrak{C}\Gamma_{m}[a]
  \mathfrak{C}\Gamma_{m}[b] \mathfrak{C}\Gamma_{m}[c]|\mathbf{I}_{m} - \mathbf{U}_{1}|
  ^{a + m} |\mathbf{I}_{m} - \mathbf{U}_{2}|^{b +m}}
  |\mathbf{C}|^{a+b+c-m} \hspace{2cm}
$$
$$
  \times \etr\left[-(\mathbf{I}_{m} - \mathbf{U}_{2})^{-1}(\mathbf{I}_{m} - \mathbf{U}_{1}\mathbf{U}_{2})
  (\mathbf{I}_{m} - \mathbf{U}_{1})^{-1}\mathbf{C}\right]
  (d\mathbf{C}) (d\mathbf{U}_{1})(d\mathbf{U}_{2}).
$$
Integrating with respect to $\mathbf{C}$ using
\begin{eqnarray*}
  \int_{\mathbf{C}=\mathbf{C}^{H} > \mathbf{0}} |\mathbf{C}|^{a+b+c-m}
  \etr\left[-(\mathbf{I}_{m} - \mathbf{U}_{2})^{-1}(\mathbf{I}_{m} - \mathbf{U}_{1}\mathbf{U}_{2})
  (\mathbf{I}_{m} - \mathbf{U}_{1})^{-1}\mathbf{C}\right] (d\mathbf{C}) \\
    = \mathfrak{C}\Gamma[a+b+c]\frac{|\mathbf{I}_{m} - \mathbf{U}_{1}|
  ^{a + b + c} |\mathbf{I}_{m} - \mathbf{U}_{2}|^{a + b+ c}}{|\mathbf{I}_{m} - \mathbf{U}_{1}\mathbf{U}_{2}|^{a + b + c
  }},
\end{eqnarray*}
(from (\ref{gamma})) gives the stated marginal density function for $
(\mathbf{U}_{1}\vdots \mathbf{U}_{2})'$. \qed

As in the real case \citep{dggj:09a}, the joint density (\ref{density1}) can be
represented as a mixture. Let us first note that
\begin{eqnarray*}
     |\mathbf{I}_{m} - \mathbf{U}_{1}\mathbf{U}_{2}|^{-(a + b + c)} &=&  {}_{1}\widetilde{F}_{0}(a + b + c
     ;\mathbf{U}_{1}\mathbf{U}_{2})\\
     &=& \sum_{t=0}^{\infty} \sum_{\tau} [a + b + c]_{\tau}
     \frac{\widetilde{C}_{\tau}(\mathbf{U}_{1}\mathbf{U}_{2})}{t!}.
\end{eqnarray*}
By substituting in (\ref{density1}) it is obtained that the joint density function of
$(\mathbf{U}_{1}\vdots \mathbf{U}_{2})'$ is
$$
   \sum_{t=0}^{\infty} \sum_{\tau} \frac{[a + b + c]_{\tau}}{\mathfrak{C}\boldsymbol{\beta}^{*}_{m}[a,b,c]} |\mathbf{U}_{1}|^{a-m}
   |\mathbf{U}_{2}|^{b-m} |\mathbf{I}_{m} - \mathbf{U}_{1}|
  |^{b + c -m} \phantom{XXXXXXXXXXXXX}
$$
\vspace{-1cm}
\begin{equation}\label{eq2}
\hspace{6cm}
  \times |\mathbf{I}_{m} - \mathbf{U}_{2}|^{a + c -m}
     \frac{\widetilde{C}_{\tau}(\mathbf{U}_{1}\mathbf{U}_{2})}{t!}.
\end{equation}
Moreover
$$
   \sum_{t=0}^{\infty} \sum_{\tau} \frac{[a + b + c]_{\tau}\mathfrak{C}\Gamma_{m}[b+c]\mathfrak{C}\Gamma_{m}[a+c]}{\mathfrak{C}\Gamma_{m}[a+b+c]
   \mathfrak{C}\Gamma_{m}[c]} \mathfrak{C}\mathcal{B}I_{m}(\mathbf{U}_{1};a,b+c) \phantom{XXXXXXXXXXXXXXXXXXXXXX}
$$
\vspace{-1cm}
$$\hspace{8cm}
  \times \  \mathfrak{C}\mathcal{B}I_{m}(\mathbf{U}_{2};b,a+c)\frac{\widetilde{C}_{\tau}(\mathbf{U}_{1}\mathbf{U}_{2})}{t!}.
$$

\section{Bimatrix variate generalised beta type II distribution}\label{sec4}

Let $\mathbf{A}$, $\mathbf{B}$ and $\mathbf{C} \in \mathfrak{C}^{m \times m}$ be
independent, where $\mathbf{A} \sim \mathfrak{C}\mathcal{G}_{m}(a, \mathbf{I}_{m})$,
$\mathbf{B} \sim \mathfrak{C}\mathcal{G}_{m}(b, \mathbf{I}_{m})$ and $\mathbf{C} \sim
\mathfrak{C}\mathcal{G}_{m}(c, \mathbf{I}_{m})$ with $\R(a)
> m-1$, $\R(b) > m-1$ and $\R(c) > m-1$ and let us define
\begin{equation}\label{bgb2}
    \mathbf{F}_{1} = \mathbf{C}^{-1/2}\mathbf{A}\mathbf{C}^{-1/2}
  \quad \mbox{and} \quad \mathbf{F}_{2} = \mathbf{C}^{-1/2}\mathbf{B}
  \mathbf{C}^{-1/2}.
\end{equation}
Clearly, $\mathbf{F}_{1} \sim \mathfrak{C}\mathcal{B}II_{m}(a,c)$ and $\mathbf{F}_{2}
\sim \mathfrak{C}\mathcal{B}II_{m}(b,c)$. But they are correlated and so the
distribution of $\mathbb{F} = (\mathbf{F}_{1}\vdots \mathbf{F}_{2})' \in \Re^{2m
\times m}$ can be termed a complex bimatrix variate generalised beta type II
distribution, which is denoted as $\mathbb{F} \sim \mathfrak{C}\mathcal{BGB}II_{2m
\times m}(a,b,c)$.

\begin{thm}\label{teob2}
Assume that $\mathbb{F} \sim \mathcal{BGB}I_{2m \times m}(a,b,c)$. Then its density
function is
\begin{equation}\label{density2}
    \frac{|\mathbf{F}_{1}|^{a-m} |\mathbf{F}_{2}|^{b-m} }{\mathfrak{C}\boldsymbol{\beta}^{*}_{m}[a,b,c]
  |\mathbf{I}_{m} + \mathbf{F}_{1} + \mathbf{F}_{2}|^{a + b + c }}(d\mathbb{F}),
\end{equation}
$\mathbf{F}_{1}= \mathbf{F}_{1}^{H} > \mathbf{0}$, $\mathbf{F}_{2} =
\mathbf{F}_{2}^{H} > \mathbf{0}$, where the measure
$$
  (d\mathbb{F}) = (d\mathbf{F}_{1})\wedge (d\mathbf{F}_{2}).
$$
and $\R(a)> m-1$, $\R(b) > m-1$ and $\R(c) > m-1$.
\end{thm}
\proof As an alternative to proceeding as in Theorem \ref{teob1}, let us recall that
if $\mathbf{U} \sim \mathfrak{C}\mathcal{B}I_{m}(a,b)$, then
$(\mathbf{I}_{m}-\mathbf{U})^{-1} -\mathbf{I}_{m} \sim
\mathfrak{C}\mathcal{B}II_{m}(a,b)$, see \citet{sk:79} and \citet{dggj:07}. Then
$$
  \mathbb{F} =
  \left (
    \begin{array}{c}
      \mathbf{F}_{1} \\
      \mathbf{F}_{2}
    \end{array}
  \right ) =
  \left (
    \begin{array}{c}
      (\mathbf{I}_{m}-\mathbf{U}_{1})^{-1}-\mathbf{I}_{m} \\
      (\mathbf{I}_{m}-\mathbf{U}_{2})^{-1}-\mathbf{I}_{m}
    \end{array}
  \right ),
$$
with the Jacobian given by (see \citet[Theorem 3.8, p. 190]{m:97})
$$
  (d\mathbf{U}_{1})(d\mathbf{U}_{2}) =  |\mathbf{I}_{m} + \mathbf{F}_{1}|^{-2m}
  |\mathbf{I}_{m} + \mathbf{F}_{2}|^{-2m} (d\mathbf{F}_{1})(d\mathbf{F}_{2}).
$$
Also, note that
\begin{eqnarray*}
  \mathbf{I}_{m}- (\mathbf{I}_{m}+\mathbf{F}_{1})^{-1}  &=& (\mathbf{I}_{m}+\mathbf{F}_{1})^{-1}
  ((\mathbf{I}_{m}+\mathbf{F}_{1}) - \mathbf{I}_{m}) = (\mathbf{I}_{m}+\mathbf{F}_{1})^{-1}\mathbf{F}_{1}\\
  \mathbf{I}_{m}- (\mathbf{I}_{m}+\mathbf{F}_{2})^{-1} &=&
  (\mathbf{I}_{m}+\mathbf{F}_{2})^{-1}\mathbf{F}_{2},
\end{eqnarray*}
Then the joint density of $(\mathbf{F}_{1}\vdots \mathbf{F}_{2})'$ is
$$
  \frac{|\mathbf{F}_{1}|^{a-m} |\mathbf{F}_{2}|^{b-m} |\mathbf{I}_{m}+\mathbf{F}_{1}|^{-(a+b+c)}
  |\mathbf{I}_{m}+\mathbf{F}_{2}|^{-(a+b+c)}}{\mathfrak{C}\boldsymbol{\beta}^{*}_{m}[a,b,c] |\mathbf{I}_{m} - (\mathbf{I}_{m} +
  \mathbf{F}_{1})^{-1}\mathbf{F}_{1} \mathbf{F}_{2}(\mathbf{I}_{m}+\mathbf{F}_{2})^{-1}|^{a + b + c
  }}(d\mathbb{F}).
$$
The desired results then follows, noting that
$$
  \frac{|\mathbf{I}_{m}+\mathbf{F}_{1}|^{-1} |\mathbf{I}_{m}+\mathbf{F}_{2}|^{-1}}
  {|\mathbf{I}_{m} - (\mathbf{I}_{m} + \mathbf{F}_{1})^{-1}\mathbf{F}_{1} \mathbf{F}_{2}
  (\mathbf{I}_{m}+\mathbf{F}_{2})^{-1}|} = |\mathbf{I}_{m} + \mathbf{F}_{1} +
  \mathbf{F}_{2}|^{-1}.   \qquad\qquad\mbox{\qed}
$$

Other properties of the distribution $\mathfrak{C}\mathcal{BGB}II_{2m \times
m}(a,b,c)$ can be found in a similar way.

\section{Properties}\label{sec5}

In this section, several basic properties and eigenvalue distributions are studied.

\subsection{Basic properties}

As direct consequences of Lemma \ref{lemma} the following basic properties are
derived: the moments $E(|\mathbf{U}_{1}|^{r}|\mathbf{U}_{2}|^{s})$, the distributions
of the product $\mathbf{Z} = \mathbf{U}_{2}^{1/2} \mathbf{U}_{1}
\mathbf{U}_{2}^{1/2}$ and the inverse $(\mathbf{U}_{1}^{-1}\vdots
\mathbf{U}_{2}^{-1})$. Their proofs are similar to those given for the real case, see
\citet{dggj:09a}.

\begin{thm}
Assume that $(\mathbf{U}_{1}\vdots\mathbf{U}_{2}) \sim \mathfrak{C}\mathcal{BGB}I_{2m
\times m}(a,b,c)$ then
$$
  E(|\mathbf{U}_{1}|^{r}|\mathbf{U}_{2}|^{s}) = \frac{\mathfrak{C}\boldsymbol{\beta}_{m}[a+r,b+c] \mathfrak{C}\boldsymbol{\beta}_{m}[b+s,a+c]}
  {\mathfrak{C}\boldsymbol{\beta}_{m}^{*}[a,b,c]}
  \phantom{XXXXXXXXXXXXXXXXXXXX}
$$
\vspace{-.5cm}
$$\phantom{XXXXXXX}
  \times {}_{3}\widetilde{F}_{2}(a+r,b+s,a+b+c;a+b+c+r,a+b+c+s;\mathbf{I}_{m}),
$$
with $\R(b+r) > m-1$, and $\R(a+c) > m-1$.
\end{thm}

\begin{thm}
Consider that $(\mathbf{U}_{1}\vdots\mathbf{U}_{2}) \sim
\mathfrak{C}\mathcal{BGB}I_{2m \times m}(a,b,c)$. Then  the density function of
$\mathbf{Z} = \mathbf{Z}^{H} = \mathbf{U}_{2}^{1/2} \mathbf{U}_{1}
\mathbf{U}_{2}^{1/2} \in \mathfrak{C}^{m \times m}$ is
$$
  \frac{\mathfrak{C}\boldsymbol{\beta}_{m}[a+c,b+c] |\mathbf{Z}|^{a-m} |\mathbf{I}_{m}-\mathbf{Z}|^{c-m}}
  {\mathfrak{C}\boldsymbol{\beta}_{m}^{*}[a,b,c]} {}_{2}\widetilde{F}_{1}(a+c,a+c;a+b+2c;\mathbf{I}_{m}-\mathbf{Z})
  (d\mathbf{Z}),
$$
and
$$
  E(|\mathbf{Z}|^{r}) = \frac{\mathfrak{C}\boldsymbol{\beta}_{m}[a+c,b+c] \boldsymbol{\beta}_{m}[a+r,c]}
  {\mathfrak{C}\boldsymbol{\beta}_{m}^{*}[a,b,c]}
  {}_{3}\widetilde{F}_{2}(c,a+c,a+c;a+c+r,a+b+2c;\mathbf{I}_{m}),
$$
with $\mathbf{0} < \R(\mathbf{Z}) < \mathbf{I}_{m}$, $\R(a+b) > m-1$ and $\R(b+c)
> m-1$.
\end{thm}

\begin{thm}
Let $(\mathbf{U}_{1}\vdots\mathbf{U}_{2})' \sim \mathfrak{C}\mathcal{BGB}I_{2m \times
m}(a,b,c)$. Then  the density function of $\mathbb{V} = (\mathbf{V}_{1}\vdots
\mathbf{V}_{2})'= ( \mathbf{U}_{1}^{-1}\vdots \mathbf{U}_{2}^{-1})' \in
\mathfrak{C}^{2m \times m}$ is
$$
   \frac{|\mathbf{V}_{1}|^{-a-m} |\mathbf{V}_{2}|^{-b-m} |\mathbf{I}_{m} - \mathbf{V}_{1}^{-1}|
  |^{b + c -m} |\mathbf{I}_{m} - \mathbf{V}_{2}^{-1}|^{a + c -m}}{\mathfrak{C}\boldsymbol{\beta}^{*}_{m}[a,b,c]
  |\mathbf{I}_{m} - (\mathbf{V}_{1}\mathbf{V}_{2})^{-1}|^{a + b + c }}(d\mathbb{V}),
$$
$\mathbf{0} < \mathbf{V}_{1} = \mathbf{V}_{1}^{H}< \mathbf{I}_{m}$, $\mathbf{0} <
\mathbf{V}_{2} = \mathbf{V}_{2}^{H}< \mathbf{I}_{m}$, where the measure
$$
  (d\mathbb{V}) = (d\mathbf{V}_{1})\wedge (d\mathbf{V}_{2}).
$$
and $\R(a)> m-1$, $\R(b) > m-1$ and $\R(c) > m-1$.
\end{thm}

\subsection{Joint eigenvalue distribution}

Many statistics of the multivariate test hypothesis are functions of the eigenvalues
or of the maximum eigenvalue. In these final two subsections, we find the eigenvalues
distribution. The following result is needed, see \citet{j:64}.

\begin{lem}\label{lem2}
If $ \ f_{_{\mathbf{A}}}(\mathbf{A}) \ (d\mathbf{A})$ is the density function of a
Hermitian matrix variate $\mathbf{A} \in \mathfrak{C}^{m \times m}$, then the
distribution of the diagonal matrix $\mathbf{\Lambda} = \diag(\lambda_{1},
\dots,\lambda_{m})$, $\lambda_{1}> \cdots > \lambda_{m} > 0$, of the eigenvalues of
$\mathbf{A}$, where $\mathbf{A} = \mathbf{G\Lambda}\mathbf{G}^{H}$, $\mathbf{G} \in
\mathcal{U}(m)$, is the eigendecomposition of $\mathbf{A}$, is
$$
  f_{_{\mathbf{\Lambda}}}(\mathbf{\Lambda}) =
  \frac{\pi^{m(m-1)}}{\mathfrak{C}\Gamma_{m}[m]}\prod_{r<s}(\lambda_{r}-\lambda_{s})^{2}
  \int_{\mathcal{U}(m)}f_{_{\mathbf{A}}}(\mathbf{G\Lambda}\mathbf{G}^{H})(d\mathbf{G}),
$$
where $(d\mathbf{G})$ is the invariant measure on the unitary group $\mathcal{U}(m)$
normalised, given as
\begin{equation}\label{unitary}
  (d\mathbf{G}) =
  \frac{\mathfrak{C}\Gamma_{m}[m]}{2^{m}\pi^{m^{2}}}(\mathbf{G}d\mathbf{G}^{H}),
  \quad \mbox{such that} \quad \int_{\mathcal{U}(m)} (d\mathbf{G}) = 1.
\end{equation}
\end{lem}

\begin{thm}
Assume that $(\mathbf{U}_{1}\vdots\mathbf{U}_{2}) \sim \mathfrak{C}\mathcal{BGB}I_{2m
\times m}(a,b,c)$ and let
\begin{equation}\label{eq3}
  \mathbb{U} =
  \left(
  \begin{array}{c}
    \mathbf{U}_{1} \\
    \mathbf{U}_{2}
  \end{array}
  \right ) =
  \left(
  \begin{array}{c}
    \mathbf{ED}_{\lambda}\mathbf{E}^{H}\\
    \mathbf{GD}_{\delta}\mathbf{G}^{H}
  \end{array}
  \right ).
\end{equation}
The spectral decomposition of $\mathbf{U}_{1}$ and $\mathbf{U}_{2}$, with
$\mathbf{E}, \mathbf{G} \in \mathcal{U}(m)$ and $\mathbf{D}_{\lambda} =
\diag(\lambda_{1}, \dots, \lambda_{m})$, $1 > \lambda_{1}> \dots, \lambda_{m} > 0$
and $\mathbf{D}_{\delta} = \diag(\delta_{1}, \dots, \delta_{m})$, $1 > \delta_{1}>
\dots, \delta_{m} > 0$. Then the joint density function of $\lambda_{1}, \dots,
\lambda_{m}, \delta_{1}, \dots, \delta_{m}$ is
$$
  \frac{\pi^{2m(m-1)}}{(\mathfrak{C}\Gamma_{m}[m])^{2}  \mathfrak{C}\boldsymbol{\beta}^{*}[a,b,c]}
  \prod_{r = 1}^{m} \left ( \lambda_{r}^{a-m} (1-\lambda_{r})^{b+c-m}\right )
  \prod_{e = 1}^{m} \left ( \delta_{e}^{b-m} (1-\delta_{e})^{a+c-m}\right )
$$
$$\times \
  \prod_{r<s}(\lambda_{r}-\lambda_{s})^{2} \prod_{e<f}(\delta_{e}-\delta_{f})^{2}
  \sum_{k = 0}^{\infty}\sum_{\kappa}\frac{[a+b+c]_{\kappa}}{k!}\frac{\widetilde{C}_{\kappa}(\mathbf{D}_{\lambda})
   \widetilde{C}_{\kappa}(\mathbf{D}_{\delta})}{\widetilde{C}_{\kappa}(\mathbf{I}_{m})}.
$$
with $\R(a)> m-1$, $\R(b) > m-1$ and $\R(c) > m-1$.
\end{thm}
\proof. The desired result follows by making the transformation (\ref{eq3}), from
Lemma \ref{lem2}, equation (\ref{unitary}) and observing that, see \citet{ch:04}
$$
  \int_{\mathbf{G} \in \mathcal{U}(m)}\int_{\mathbf{E} \in \mathcal{U}(m)}
  \widetilde{C}_{\kappa} (\mathbf{ED}_{\lambda}\mathbf{E}^{H}
  \mathbf{GD}_{\delta}\mathbf{G}^{H}) (d\mathbf{E})(d\mathbf{G})= \frac{
  \widetilde{C}_{\kappa} (\mathbf{D}_{\lambda})\widetilde{C}_{\kappa}
  (\mathbf{D}_{\delta})}{\widetilde{C}_{\kappa} (\mathbf{I})}. \mbox{\qed}
$$

\subsection{Joint distribution of $\lambda_{max}$ and $\delta_{max}$}

In this subsection, we derive the distribution of the largest eigenvalues,
$\lambda_{max}$ and $\delta_{max}$ of a complex bimatrix variate beta type I matrix.
With this aim, let us consider the following theorem.

\begin{thm}\label{teo5}
Assume that $(\mathbf{U}_{1}\vdots\mathbf{U}_{2}) \sim \mathfrak{C}\mathcal{BGB}I_{2m
\times m}(a,b,c)$ and let $\mathbf{\Delta}_{1}, \mathbf{\Delta}_{2} \in
\mathfrak{C}^{m \times m}$ be Hermitian positive definite matrices, $\mathbf{0} <
\mathbf{\Delta}_{1}, \mathbf{\Delta}_{2} < \mathbf{I}$. Then the probability
$P(\mathbf{U}_{1} < \mathbf{\Delta}_{1}, \mathbf{U}_{2} < \mathbf{\Delta}_{2})$ is
given by
$$
  \frac{\mathfrak{C}\boldsymbol{\beta}_{m}[a,m]
  \mathfrak{C}\boldsymbol{\beta}_{m}[b,m]}
  {\mathfrak{C}\boldsymbol{\beta}_{m}^{*}[a,b,c]} |\mathbf{\Delta}_{1}|^{a}
  |\mathbf{\Delta}_{2}|^{b} \sum_{k,t,s =0}^{\infty} \  \sum_{\kappa,\tau,\sigma; \phi_{1} \in
  \kappa.\tau; \phi_{2} \in \sigma^{*}.\phi_{1}^{*}} \frac{[a]_{\phi_{2}}[b]_{\phi_{1}}
  [a+b+c]_{\kappa}}{k! \ t! \ s!}
$$
$$\hspace{1.5cm}\times \
  \frac{[-(a+c)+m]_{\tau} [-(b+c)+m]_{\sigma}} {[a+m]_{\phi_{2}} [b+m]_{\phi_{1}}}
  \theta_{\phi_{1}}^{\kappa, \tau} \pi_{\sigma,\phi_{1}}^{\sigma, \tau, \kappa;
  \phi_{2}}\widetilde{C}_{\phi_{2}}^{\sigma, \tau, \kappa}(\mathbf{\Delta}_{1},
  \mathbf{\Delta}_{2},\mathbf{\Delta}_{1} \mathbf{\Delta}_{2})
$$
where
$$
  \sum_{\kappa,\tau,\sigma; \phi_{1} \in  \kappa.\tau; \phi_{2} \in
  \sigma^{*}.\phi_{1}^{*}} = \sum_{\kappa,\tau,\sigma} \ \ \sum_{\phi_{1} \in
  \kappa.\tau} \ \sum_{\phi_{2} \in \sigma^{*}.\phi_{1}^{*}}
$$
and $\R(a)> m-1$, $\R(b) > m-1$ and $\R(c) > m-1$.
\end{thm}
\proof From (\ref{density1}) the probability $P(\mathbf{U}_{1} < \mathbf{\Delta}_{1},
\mathbf{U}_{2} < \mathbf{\Delta}_{2})$ is
\begin{equation}\label{eq4}
    \int_{\mathbf{O}}^{\mathbf{\Delta}_{1}} \int_{\mathbf{O}}^{\mathbf{\Delta}_{2}}
  \frac{|\mathbf{U}_{1}|^{a-m} |\mathbf{U}_{2}|^{b-m} |\mathbf{I}_{m} - \mathbf{U}_{1}|
  |^{b + c -m} |\mathbf{I}_{m} - \mathbf{U}_{2}|^{a + c -m}}{\mathfrak{C}\boldsymbol{\beta}^{*}_{m}[a,b,c]
  |\mathbf{I}_{m} - \mathbf{U}_{1}\mathbf{U}_{2}|^{a + b + c
  }}(d\mathbf{U}_{1})(d\mathbf{U}_{}).
\end{equation}
By rewritten (\ref{eq4}) as in (\ref{eq2}), and integrating with respect to
$\mathbf{U}_{2}$ using (\ref{ib1}), it is obtained that
$$
  \frac{\mathfrak{C}\boldsymbol{\beta}_{m}[b,m]}
  {\mathfrak{C}\boldsymbol{\beta}_{m}^{*}[a,b,c]} |\mathbf{\Delta}_{1}|^{a}
  \sum_{k,t =0}^{\infty} \  \sum_{\kappa,\tau; \phi_{1} \in
  \kappa.\tau} \frac{[b]_{\phi_{1}}[a+b+c]_{\kappa}[-(a+c)+m]_{\tau}\theta_{\phi_{1}}^{\kappa, \tau}}{k! \ t! [b+m]_{\phi_{1}}}
$$
$$\hspace{1.5cm}\times \
  \int_{\mathbf{O}}^{\mathbf{\Delta}_{1}} |\mathbf{U}_{1}|^{a-m} |\mathbf{I}_{m} -
  \mathbf{U}_{1}| |^{b + c -m} \widetilde{C}_{\phi_{1}}^{\tau, \kappa}
  (\mathbf{\Delta}_{2},\mathbf{\Delta}_{2} \mathbf{U}_{1})(d\mathbf{U}_{1}).
$$
The desired result follows by integrating with respect to $\mathbf{U}_{1}$ using
(\ref{ib2}). \qed

The following result is obtained from Theorem \ref{teo5}.
\begin{cor}
Let $(\mathbf{U}_{1}\vdots\mathbf{U}_{2}) \sim \mathfrak{C}\mathcal{BGB}I_{2m \times
m}(a,b,c)$, $\R(a)> m-1$, $\R(b) > m-1$ and $\R(c) > m-1$. If $\lambda_{max}$ and
$\delta_{max}$ are the largest eigenvalues of $(\mathbf{U}_{1}$ and $\mathbf{U}_{2}$,
respectively, then their joint distribution function, $P(\lambda_{max} < x,
\delta_{max} <y)$ is given by
$$
  \frac{\mathfrak{C}\boldsymbol{\beta}_{m}[a,m]
  \mathfrak{C}\boldsymbol{\beta}_{m}[b,m]}
  {\mathfrak{C}\boldsymbol{\beta}_{m}^{*}[a,b,c]} \sum_{k,t,s =0}^{\infty} \  \sum_{\kappa,\tau,\sigma; \phi_{1} \in
  \kappa.\tau; \phi_{2} \in \sigma^{*}.\phi_{1}^{*}} \frac{x^{am+s+k}  y^{bm+t+k}[a]_{\phi_{2}}[b]_{\phi_{1}}
  }{k! \ t! \ s!}
$$
$$\hspace{1.5cm}\times \
  \frac{[-(a+c)+m]_{\tau} [-(b+c)+m]_{\sigma}} {[a+b+c]_{\kappa}[a+m]_{\phi_{2}} [b+m]_{\phi_{1}}}
  \theta_{\phi_{1}}^{\kappa, \tau} \pi_{\sigma,\phi_{1}}^{\sigma, \tau, \kappa;
  \phi_{2}} \theta_{\phi_{2}}^{\kappa, \tau, \sigma} \widetilde{C}_{\phi_{2}}^{\sigma, \tau,
  \kappa}(\mathbf{I}).
$$
\end{cor}
\proof. The inequalities $\lambda_{max} < x$ and  $\delta_{max} <y$ are equivalent to
$\mathbf{U}_{1} < x\mathbf{I}$ and  $\mathbf{U}_{2} < y\mathbf{I}$. Therefore, the
result follows by letting $\mathbf{\Delta}_{1} = x\mathbf{I}$ and
$\mathbf{\Delta}_{2} = y\mathbf{I}$ in Theorem \ref{teo5}, and taking into account
that from \citet[eqs. (2.1) and (2.7)]{da:79}
$$
  \widetilde{C}_{\phi_{2}}^{\sigma, \tau, \kappa}(x\mathbf{I},
  y \mathbf{I},xy\mathbf{I}) = x^{s+k}  y^{t+k}\widetilde{C}_{\phi_{2}}^{\sigma, \tau,
  \kappa}(\mathbf{I}). \qquad\qquad \mbox{\qed}
$$

\section{Conclusions}

Random matrices play a prominent role because of their deep mathematical structure.
They have arisen in a number of fields (statistics, graph theory, stochastic linear
algebra, physics, signal processing, etc.), often independently. In particular, the
statisticians are interested in studying the Normal (Gaussian), Wishart,
\textsc{manova}, and circular random matrices, which, from the point of view of
\emph{random matrix theory} are termed Hermite, Laguerre, Jacobi, and Fourier
ensembles, \citet{m:91} and \citet{e:05}.

The results obtained in the present work can be considered as a generalisation of the
Jacobi ensemble for the case in which there are two correlated Jacobi ensembles, and
with all their potential applications, see \citet{es:08}.

Another potential use appears in the context of complex shape theory, see
\citet{mdm:06}. Specifically, in the approach known as affine shape or configuration
densities, see \citet{cdggf:09}.

\section*{Acknowledgments}

This research work was partially supported  by CONACYT-M\'exico, Research Grant No. \
81512 and IDI-Spain, Grants No. FQM2006-2271 and MTM2008-05785. This paper was
written during J. A. D\'{\i}az- Garc\'{\i}a's stay as a visiting professor at the
Department of Statistics and O. R. of the University of Granada, Spain.

\end{document}